\documentclass[11pt]{article}
\usepackage{latexsym,amsmath,amsfonts,graphicx,pmat}
\begin{document}
\begin{center}
{\LARGE\textbf{More on the normalized Laplacian Estrada index}}\\
\bigskip
\bigskip
Yilun Shang\\
Singapore University of Technology and Design\\
20 Dover Drive, Singapore 138682\\
e-mail: \texttt{shylmath@hotmail.com}
\end{center}

\smallskip
\begin{abstract}

Let $G$ be a simple graph of order $N$. The normalized Laplacian
Estrada index of $G$ is defined as
$NEE(G)=\sum_{i=1}^Ne^{\lambda_i}$, where
$\lambda_1,\lambda_2,\cdots,\lambda_N$ are the normalized Laplacian
eigenvalues of $G$. In this paper, we give a tight lower bound for
$NEE$ of general graphs. We also calculate $NEE$ for a class of
treelike fractals, which contain some classical chemical trees as
special cases. It is shown that $NEE$ scales linearly with the order
of the fractal, in line with a best possible lower bound for
connected bipartite graphs.

\bigskip

\bigskip


\end{abstract}

\bigskip
\normalsize

\section{Introduction}

Let $G$ be a simple undirected graph with vertex set
$V(G)=\{v_1,v_2\cdots,v_N\}$. Denote by
$A=A(G)\in\mathbb{R}^{N\times N}$ the adjacency matrix of $G$, and
$\lambda_1(A)\ge\lambda_2(A)\ge\cdots\ge\lambda_N(A)$ the
eigenvalues of $A$ in non-increasing order. The Laplacian and
normalized Laplacian matrices of $G$ are defined as $L=L(G)=D-A$ and
$\mathcal{L}=\mathcal{L}(G)=D^{-1/2}LD^{-1/2}$, respectively, where
$D=D(G)$ is a diagonal matrix with $d_i$ on the main diagonal being
the degree of vertex $v_i$, $i=1,\cdots,N$. Here, by convention
$d_i^{-1}=0$ if $d_i=0$. The eigenvalues of $L$ and $\mathcal{L}$
are referred to as the Laplacian and normalized Laplacian
eigenvalues of graph $G$, denoted by
$\lambda_1(L)\ge\lambda_2(L)\ge\cdots\ge\lambda_N(L)$ and
$\lambda_1(\mathcal{L})\ge\lambda_2(\mathcal{L})\ge\cdots\ge\lambda_N(\mathcal{L})$,
respectively. Details of the theory of graph eigenvalues can be
found in \cite{1,2}.

The Estrada index of a graph $G$ is defined in \cite{3} as
$$
EE=EE(G)=\sum_{i=1}^Ne^{\lambda_i(A)},
$$
which was first introduced in 2000 as a molecular
structure-descriptor by Ernesto Estrada \cite{4}. The Laplacian
Estrada index of a graph $G$ is defined in \cite{6} as
$$
LEE=LEE(G)=\sum_{i=1}^Ne^{(\lambda_i(L)-2E/N)},
$$
where $E$ is the number of edges in $G$. Another essentially
equivalent definition for the Laplacian Estrada index is given by
$LEE=\sum_{i=1}^Ne^{\lambda_i(L)}$ in \cite{7} independently. Albeit
young, the (Laplacian) Estrada index has already found a variety of
chemical applications in the degree of folding of long-chain
polymeric molecules \cite{5,8}, extended atomic branching \cite{9},
and the Shannon entropy descriptor \cite{10}. In addition,
noteworthy applications in complex networks are uncovered
\cite{11,12,13,14}. Mathematical properties of $EE$ and $LEE$,
especially the upper and lower bounds, are investigated in e.g.
\cite{15,16,17,18,19,20,21,27}.

Very recently, in full analogy with the (Laplacian) Estrada index,
the normalized Laplacian Estrada index of a graph $G$ is introduced
in \cite{22} as
\begin{equation}
NEE=NEE(G)=\sum_{i=1}^Ne^{\lambda_i(\mathcal{L})-1}.\label{1}
\end{equation}
Among other things, the following tight bounds for $NEE$ are
obtained.

\smallskip
\noindent\textbf{Theorem 1.} \cite{22} \itshape \quad If $G$ is a
connected graph of order $N$, then
\begin{equation}
NEE\ge(N-1)e^{1/(N-1)}+e^{-1}.\label{2}
\end{equation}
The equality holds if and only if $G=K_N$, i.e., a complete graph.
\normalfont
\smallskip

\noindent\textbf{Theorem 2.} \cite{22} \itshape \quad Let $G$ be a
connected bipartite graph of order $N$ with maximum degree $\Delta$
and minimum degree $\delta$. Then
\begin{equation}
e^{-1}+e+\sqrt{(N-2)^2+\frac{2(N-2\Delta)}{\Delta}}\le NEE\le
e^{-1}+e+(N-3)-\sqrt{\frac{N-2\delta}{\delta}}+e^{\sqrt{\frac{N-2\delta}{\delta}}}.\label{3}
\end{equation}
The equality occurs in both bounds if and only if $G$ is a complete
bipartite regular graph. \normalfont
\smallskip

In this paper, we find a tight lower bound for $NEE$ of a general
graph (not necessarily connected) by extending Theorem 1. We also
calculate $NEE$ for a class of treelike fractals, which subsume some
important chemical trees as special cases, through an explicit
recursive relation. We unveil that $NEE$ scales linearly with the
order of the fractal, i.e., $NEE\propto N$ for large $N$, in line
with the lower bound in Theorem 2.

\section{A new lower bound for general graphs}

We begin with some basic properties of the normalized Laplacian
eigenvalues of a connected graph $G$. In the rest of the paper, to
ease notation, we shall use $\lambda_i$, $i=1,\cdots,N$ to represent
the normalized Laplacian eigenvalues of a graph of order $N$.

\smallskip
\noindent\textbf{Lemma 1.} \cite{2} \itshape \quad Let $G$ be a
connected graph of order $N\ge2$. Then
\begin{itemize}
\item[(i)] $\sum_{i=1}^N\lambda_i=N$.
\item[(ii)] $\lambda_N=0$ and $\lambda_i\in(0,2]$ for every
$1\le i\le N-1$.
\item[(iii)] If $G$ is a bipartite graph, then $\lambda_1=2$ and
$\lambda_2<2$.
\end{itemize}
\normalfont
\smallskip

\noindent\textbf{Theorem 3.} \itshape \quad Let $G$ be a graph of
order $N$. If $G$ possesses $c$ connected components, $r$ of which
are isolated vertices, then
\begin{equation}
NEE\ge(N-c)e^{\frac{c-r}{N-c}}+ce^{-1}.\label{4}
\end{equation}
The equality holds if and only if $G$ is a union of copies of $K_s$,
for some fixed $s\ge2$, and $r$ isolated vertices. \normalfont
\smallskip

Clearly, we reproduce Theorem 1 when $c=1$ (and $r=0$).

\noindent\textbf{Proof.} From Lemma 1 and the fact that the
normalized Laplacian eigenvalue for an isolated vertex is zero, we
obtain $\lambda_N=\cdots=\lambda_{N-c+1}=0$ and
$\lambda_1+\cdots+\lambda_{N-c}=N-r$. Hence,
$$
NEE=e^{-1}\left(c+\sum_{i=1}^{N-c}e^{\lambda_i}\right)\ge
e^{-1}\left(c+(N-c)e^{\frac{\lambda_1+\cdots+\lambda_{N-c}}{N-c}}\right)=ce^{-1}+(N-c)e^{\frac{c-r}{N-c}},
$$
where the arithmetic-geometric mean inequality is utilized.

Now suppose $G=(c-r)K_s\cup rK_1$ for some integers $0\le r\le c$
and $s\ge2$. Then $N=(c-r)s+r=cs-r(s-1)$, and the normalized
Laplacian eigenvalues of $G$ can be computed as 0 with multiplicity
$c$ and $s/(s-1)$ with multiplicity $(c-r)(s-1)$. Therefore,
$$
NEE(G)=ce^{-1}+(c-r)(s-1)e^{\frac{s}{s-1}-1}=ce^{-1}+(N-c)e^{\frac{c-r}{N-c}}.
$$

Conversely, suppose that the equality holds in (\ref{4}). Then from
the above application of the arithmetic-geometric mean inequality we
know that all the non-zero normalized Laplacian eigenvalues of $G$
must be mutually equal. Suppose that a connected graph $H$ with
order $s$ has normalized Laplacian eigenvalues $\lambda>0$ with
multiplicity $s-1$ and a single zero eigenvalues. It now suffices to
show $H=K_s$.

The case of $s=1$ holds trivially. In what follows, we assume
$s\ge2$. We first claim that
\begin{equation}
\mathcal{L}(H)(\mathcal{L}(H)-\lambda I_s)={\bf0},\label{5}
\end{equation}
where $I_s\in\mathbb{R}^{s\times s}$ is the identity matrix and
$\bf0$ is the zero matrix. Indeed, for any vector
$x\in\mathbb{R}^s$, we can decompose it as
$x=a_1D^{1/2}{\bf1}+\sum_{i=2}^sa_ix_i$, where $a_i\in\mathbb{R}$,
$\bf1$ is a column vector with all elements being 1, and $x_i$,
$i=2,\cdots,s$ are eigenvectors associated with $\lambda$. Thus, the
equation (\ref{5}) follows by the fact that
$\mathcal{L}(H)(\mathcal{L}(H)-\lambda I_s)x=0$ for all $x$.

In the light of (\ref{5}), all columns of $\mathcal{L}(H)-\lambda
I_s$ belong to the null space of $\mathcal{L}(H)$. It follows from
Lemma 1 that the null space of $\mathcal{L}(H)$ is of dimension one
and is spanned by $D^{1/2}\bf1$. Therefore, each column of
$\mathcal{L}(H)-\lambda I_s$ is of the form $\alpha D^{1/2}\bf1$ for
some $\alpha\in\mathbb{R}$. Since $H$ is connected, $\alpha\not=0$
and $d_i>0$ for all $i$. It follows that any pair of vertices in $H$
must be adjacent in view of the definition of $\mathcal{L}(H)$. This
means $H=K_s$. $\Box$

\section{The normalized Laplacian Estrada index of treelike fractals}

In this section, we analytically calculate $NEE$ for a class of
treelike fractals through a recursive relation deduced from the
self-similar structure of the fractals. It is of great interest to
seek $NEE$ (also $LEE$ and $EE$) for specific graphs due to the
following two reasons: Firstly, a universal approach for evaluating
these structure descriptors of general graphs, especially
large-scale graphs, is out of reach so far; and secondly, the known
upper and lower bounds (such as (\ref{3})) are too far away apart to
offer any meaningful guide for a given graph.

\begin{figure}[hbt]
\centering
\scalebox{0.45}{\includegraphics[3pt,4pt][511pt,174pt]{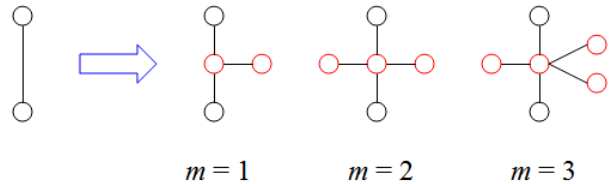}}
\caption{Building blocks of the fractals. The next generation is
obtained from current generation through replacing each edge by the
stars on the right-hand side of the arrow.}
\end{figure}

The fractal graphs in question are constructed in an iterative
manner \cite{23}. For integers $n\ge0$ and $m\ge1$, let $G_n(m)$
represent the graph after $n$ iterations (generations). When $n=0$,
$G_0(m)$ is an edge linking two vertices. In each following
iteration $n\ge1$, $G_n(m)$ is built from $G_{n-1}(m)$ by conducting
such operations on each existing edges as shown in Fig. 1: subdivide
the edge into two edges, connecting to a new vertex; then, generate
$m$ new vertices and attach each of them to the middle vertex of the
length-2 path. In Fig. 2 are illustrated the first several steps of
the iterative construction process corresponding to $m=1$. Clearly,
it reduces to two classes of chemical graphs \cite{24}: the $T$
fractal (when $m=1$) and the Peano basin fractal (when $m=2$).

\begin{figure}[hbt]
\centering
\scalebox{0.45}{\includegraphics[15pt,6pt][466pt,469pt]{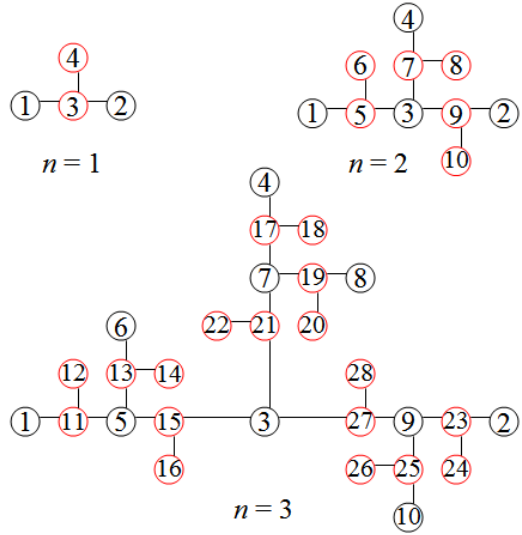}}
\caption{Growth process for a fractal corresponding to $m=1$. The
vertex with index 3 is called the \textit{inmost} vertex since it is
in the center of the graph, and those vertices farthest from the
inmost vertex are called the \textit{outmost} vertices (e.g.
vertices with indices 1, 2, and 4).}
\end{figure}

Another intuitive generation approach of the fractal $G_n(m)$, which
will be used later, highlights the self-similarity. Taking $G_n(m)$
with $n=3$ and $m=1$ as an example (see Fig. 2), $G_n(m)$ can be
obtained by coalescing $m+2$ replicas of $G_{n-1}(m)$ (denoted by
$G_{n-1}^{(i)}(m)$, $i=1,\cdots,m+2$) with the $m+2$ outmost
vertices in separate duplicates being merged into one single new
vertex --- the inmost vertex of $G_n(m)$ (e.g. the vertex with index
3 in Fig. 2).

Since our results will be stated for any fixed $m$, we often
suppress the index $m$ in notations. Some basic properties of
$G_n=G_n(m)$ are easy to derive. For example, the number of vertices
and edges are given by $N_n=(m+2)^n+1$ and $E_n=(m+2)^n$,
respectively. We write $\mathcal{L}_n=\mathcal{L}(G_n)$ as the
normalized Laplacian matrix of $G_n$. Its eigenvalues are denoted by
$\lambda_1(n)\ge\lambda_2(n)\ge\cdots\ge\lambda_{N_n}(n)$.
Therefore, the normalized Laplacian Estrada index
\begin{equation}
NEE(G_n)=\sum_{i=1}^{N_n}e^{\lambda_i(n)-1}\label{6}
\end{equation}
can be easily derived provided we have all the eigenvalues.

\smallskip
\noindent\textbf{Theorem 4.} \itshape \quad All the eigenvalues
$\{\lambda_i(n)\}_{i=1}^{N_n}$, $n\ge1$ can be obtained by the
following recursive relations:
\begin{itemize}
\item[(i)] 0 and 2 are both single eigenvalues for every $n\ge1$.

\item[(ii)] 1 is an eigenvalue with multiplicity $m(m+2)^{n-1}+1$ for $n\ge1$.

\item[(iii)] For $n\ge1$, all eigenvalues $\lambda(n+1)$ (except 0,1, and 2) at generation $n+1$
are exactly those produced via
\begin{equation}
\lambda(n+1)=1\pm\sqrt{1-\frac{\lambda(n)}{m+2}}\label{7}
\end{equation}
by using eigenvalues $\lambda(n)$ (except 0 and 2) at generation
$n$.
\end{itemize}
\normalfont
\smallskip

\noindent\textbf{Proof.} We start with checking the completeness of
the eigenvalues provided by the rules (i), (ii), and (iii). It is
direct to check that the eigenvalues for $G_1(m)$ (0, 2, and 1 with
multiplicity $m+1$) given by (i) and (ii) are complete. Therefore,
all eigenvalues (except 0,1, and 2) for $G_n(m)$, $n\ge2$ are
descendants of eigenvalue 1 following (\ref{7}). Each father
eigenvalue produces 2 child eigenvalues in the next generation.
Thus, the total number of eigenvalues of $G_n(m)$ is found to be
$$
2+\sum_{i=0}^{n-1}\left(m(m+2)^{n-1-i}+1\right)2^i=(m+2)^n+1=N_n,
$$
which implies that all eigenvalues are obtained.

Since $G_n$ is connected and bipartite, (i) follows by Lemma 1. It
suffices to show (ii) and that each eigenvalue $\lambda_i(n+1)$
(except 0,1, and 2) can be derived through (\ref{7}) by some
eigenvalue $\lambda_i(n)$. Since $\mathcal{L}_n$ is similar to
$I_{N_n}-D^{-1}(G_n)A(G_n)$ thus having the same eigenvalues, we
will focus on $I_{N_n}-D^{-1}(G_n)A(G_n)$ in the sequel and write it
as $\mathcal{L}_n$ for simplicity.

To derive the recursive relation (\ref{7}), we resort to the
so-called decimation method \cite{25}. Let $\alpha$ denote the set
of vertices belonging to $G_n$ and $\beta$ the set of vertices
created at iteration $n+1$. Assume that $\lambda_i(n+1)$ is an
eigenvalue of $\mathcal{L}_{n+1}$ and $\lambda_i(n+1)\not=0,1,2$.
Then the eigenvalue equation for $\mathcal{L}_{n+1}$ can be recast
in the following block form
\begin{equation}
\mathcal{L}_{n+1}\left[\begin{array}{c}u_{\alpha}\\
u_{\beta}\end{array}\right]=\left[\begin{array}{cc}
I_{N_n}&\mathcal{L}_{\alpha,\beta}\\\mathcal{L}_{\beta,\alpha}&\mathcal{L}_{\beta,\beta}
\end{array}\right]\left[\begin{array}{c}u_{\alpha}\\
u_{\beta}\end{array}\right]=\lambda_i(n+1)\left[\begin{array}{c}u_{\alpha}\\
u_{\beta}\end{array}\right],\label{8}
\end{equation}
where $\mathcal{L}_{\beta,\beta}=I_{E_n}\otimes B$ with
$$
B=\left[\begin{array}{ccccc}1&-\frac{1}{m+2}&-\frac{1}{m+2}&\cdots&-\frac{1}{m+2}\\
-1&1&0&\cdots&0\\-1&0&1&\cdots&0\\
\vdots&\vdots&\vdots&\ddots&\vdots\\
-1&0&0&\cdots&1\end{array}\right]\in\mathbb{R}^{(m+1)\times(m+1)}.
$$
By eliminating $u_{\beta}$ from (\ref{8}) we arrive at
\begin{equation}
\left(I_{N_n}+\mathcal{L}_{\alpha,\beta}(\lambda_i(n+1)I_{(m+1)E_n}-\mathcal{L}_{\beta,\beta})^{-1}\mathcal{L}_{\beta,\alpha}\right)u_{\alpha}=\lambda_i(n+1)u_{\alpha},\label{9}
\end{equation}
provided the concerned matrix is invertible. Let
$\mathcal{P}_n=I_{N_n}+\mathcal{L}_{\alpha,\beta}(\lambda_i(n+1)I_{(m+1)E_n}-\mathcal{L}_{\beta,\beta})^{-1}\mathcal{L}_{\beta,\alpha}$
and $\mathcal{Q}_n=(x+y)I_{N_n}-y\mathcal{L}_n$, with
\begin{equation}
x=\frac{\lambda_i(n+1)-1}{(m+2)(\lambda_i(n+1)-2)\lambda_i(n+1)+2}+1=y+1.\label{10}
\end{equation}
It is not difficult to see that $\mathcal{P}_n=\mathcal{Q}_n$.

Indeed,
\begin{eqnarray*}
(\lambda_i(n+1)I_{(m+1)E_n}-\mathcal{L}_{\beta,\beta})^{-1}&=&I_{E_n}\otimes(\lambda_i(n+1)I_{m+1}-B)^{-1}\\
&=&I_{E_n}\otimes\frac{\mathrm{adj}((\lambda_i(n+1)I_{m+1}-B))}{\det((\lambda_i(n+1)I_{m+1}-B)},
\end{eqnarray*}
where $\mathrm{adj}(\cdot)$ means the adjugate matrix. Let $z$ be
the element on the first row and the first column of
$(\lambda_i(n+1)I_{m+1}-B)^{-1}$. We have
$$
z=\frac{(\lambda_i(n+1)-1)^m}{\det((\lambda_i(n+1)I_{m+1}-B)}=\frac{(m+2)(\lambda_i(n+1)-1)}{2+(m+2)\lambda_i(n+1)(\lambda_i(n+1)-2)},
$$
which is well-defined since $\lambda_i(n+1)\not=0,2$. With these
preparations, it is easy to make an entry-wise comparison between
$\mathcal{P}_n$ and $\mathcal{Q}_n$. Clearly,
$(\mathcal{Q}_n)_{i,i}=x=1+d_i(-1/d_i)z(-1/(m+2))=(\mathcal{P}_n)_{i,i}$;
for $i\not=j$, if $i$ and $j$ are not adjacent,
$(\mathcal{Q}_n)_{i,j}=(\mathcal{P}_n)_{i,j}=0$, while if $i$ and
$j$ are adjacent,
$(\mathcal{Q}_n)_{i,j}=-y(-1/d_i)=(-1/d_i)z(-1/(m+2))=(\mathcal{P}_n)_{i,j}$.
Hence, we conclude $\mathcal{P}_n=\mathcal{Q}_n$.

Inserting the equality $\mathcal{P}_n=\mathcal{Q}_n$ into (\ref{9})
we get
$$
\mathcal{L}_nu_{\alpha}=\left(\frac{x+y-\lambda_i(n+1)}{y}\right)u_{\alpha},
$$
where $y\not=0$ since $\lambda_i(n+1)\not=1$. This indicates that
$\lambda_i(n)=(x+y-\lambda_i(n+1))/y$, where $\lambda_i(n)$ is an
eigenvalue of $\mathcal{L}_n$ associated with the eigenvector
$u_{\alpha}$. Combining this with (\ref{10}) yields a quadratic
equation, whose solution gives the formula (\ref{7}) as desired.

It remains to show (ii). Let $M_n(\lambda)$ represent the
multiplicity of eigenvalue $\lambda$ of $\mathcal{L}_n$. We have
$$
M_n(\lambda=1)=N_n-\mathrm{rank}(\mathcal{L}_n-I_{N_n}).
$$
The problem of determining multiplicity is reduced
to evaluating $\mathrm{rank}(\mathcal{L}_n-I_{N_n})$. (ii) follows
by showing that
\begin{equation}
\mathrm{rank}(\mathcal{L}_n-I_{N_n})=2(m+2)^{n-1}\label{11}
\end{equation}
for each $n\ge1$.

To show (\ref{11}) we use the method of induction. For $n=1$, we
have
$$
\mathcal{L}_1-I_{N_1}=\left[\begin{array}{cccc}
0&-\frac{1}{m+2}&\cdots&-\frac{1}{m+2}\\-1&0&\cdots&0\\
\vdots&\vdots&\ddots&\vdots\\-1&0&\cdots&0
\end{array}\right]\in\mathbb{R}^{(m+3)\times(m+3)}.
$$
Thus, $\mathrm{rank}(\mathcal{L}_1-I_{N_1})=2$. For $n=2$, we have
$$
\mathcal{L}_2-I_{N_2}=\begin{pmat}[{|.}]
0&-\frac{1}{m+2}&\cdots&-\frac{1}{m+2}\cr\- 0&&&\cr
\vdots&&I_{m+2}\otimes F&\cr 0&&&\cr
\end{pmat},
$$
where $F\in\mathbb{R}^{(m+2)\times(m+2)}$ is the matrix
$\mathcal{L}_1-I_{N_1}$ deleting the last column and the last row.
Accordingly, $\mathrm{rank}(\mathcal{L}_2-I_{N_2})=2(m+2)$. For
$n\ge2$, we have
$$
\mathcal{L}_{n+1}-I_{N_{n+1}}=\begin{pmat}[{|.}]
0&-\frac{1}{m+2}&\cdots&-\frac{1}{m+2}\cr\- 0&&&\cr
\vdots&&I_{m+2}\otimes F_n&\cr 0&&&\cr
\end{pmat},
$$
where $F_n+I_{(m+2)^n}=\mathcal{L}(G_n\backslash\{\mathrm{an}\
\mathrm{outmost}\ \mathrm{vertex}\})$. Moreover, $F_n$ can be
iteratively expressed as
$$
F_n=\left[\begin{array}{ccccc}F_{n-1}&0&\cdots&0&w_1\\
0&F_{n-1}&\cdots&0&w_2\\
\vdots&\vdots&\ddots&\vdots&\vdots\\
0&0&\cdots&F_{n-1}&w_{m+1}\\
w_1^T&w_2^T&\cdots&w_{m+1}^T&F_{n-1}\end{array}\right],
$$
where each $w_i\in\mathbb{R}^{(m+2)^{n-1}\times(m+2)^{n-1}}$ is a
matrix containing only one non-zero element $-1/(m+2)$ describing
the edge linking the inmost vertex in $G_n$ to one vertex in the
replica $G_{n-1}^{(i)}$. For any vertex $u$ in $G_n$ that is
adjacent to the inmost vertex of $G_n$, it has a neighbor $v$ with
degree one. Hence, there is only one non-zero element for row $v$
and for column $v$, respectively, that is, $(F_n)_{v,u}=-1$ and
$(F_n)_{u,v}=-1/(m+2)$. By using some basic operations for the
matrix, we can eliminate all non-zero elements at the last row and
the last column of $F_n$. Thus,
$\mathrm{rank}(\mathcal{L}_{n+1}-I_{N_{n+1}})=(m+2)\cdot\mathrm{rank}(F_n)=(m+2)\cdot\mathrm{rank}(\mathcal{L}_n-I_{N_n})$.
This yields (\ref{11}), and finally concludes the proof. $\Box$

\noindent\textbf{Remark 1.}  We mention that although the
eigenvalues of a related matrix of $G_n(m)$ have been computed in
\cite{26} by a semi-analytical method, the normalized Laplacian
eigenvalues cannot be derived directly from results therein.

With Theorem 4 at hand, the normalized Laplacian Estrada index can
be easily evaluated through (\ref{6}). Note that $G_n(m)$ is a
connected bipartite graph with maximum degree $\Delta=m+2$ and
minimum degree $\delta=1$. In Fig. 3, we display $NEE(G_n(m))$ for
$n=1,\cdots,7$, $m=1,\cdots,5$ together with the obtained bounds in
Theorem 2. The results gathered in Fig. 3 allow us to draw several
interesting comments. First, as expected from Theorem 2, all values
of $NEE(G_n(m))$ lie between the upper and lower bounds. Second, the
lower bounds for different values of $\Delta$ are collapsed together
and they scale with the order of the graph as $NEE\propto N$, which
can be derived from (\ref{3}). Similarly, the upper bound scales
with $N$ as $NEE\propto e^{\sqrt{N}}$. Third, $NEE(G_n(m))$ also
scales linearly with the order of the fractal, i.e.,
$NEE(G_n(m))\propto N_n$, in parallel with the lower bound.

\begin{figure}[htb]
\centering{\includegraphics[width=.8\textwidth]{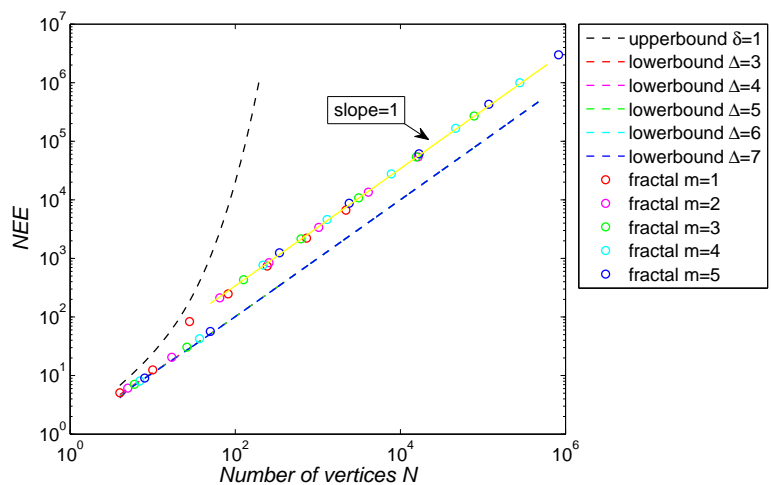}}
\caption{The normalized Laplacian Estrada index $NEE$ versus the
number of vertices $N$ in log-log scale. Circles represent
$NEE(G_n(m))$ for $n=1,\cdots,7$, $m=1,\cdots,5$; dashed lines
represent the bounds given by (\ref{3}).}
\end{figure}

\end{document}